\documentclass[12pt,a4paper]{article}
\usepackage[dvips]{graphicx}
\usepackage{amsmath, amssymb,latexsym}
\usepackage{mathrsfs}
\usepackage{psfrag}
\usepackage[all]{xy}
\usepackage{color}
\usepackage{hyperref}
\usepackage{enumitem}
\usepackage{lineno}
\usepackage[T1]{fontenc}
\usepackage{mathtools,halloweenmath}
\usepackage{textcomp}
\newtheorem{prop}{Proposition}[section]
\newtheorem{theo}[prop]{Theorem}
\newtheorem{cor}[prop]{Corollary}

\newtheorem{defn}[prop]{Definition}
\newtheorem{rem}[prop]{Remark}

\newenvironment{proof}
{\begin{trivlist} \item[\hskip \labelsep {\bf Proof}\hspace*{3 mm}]}
	{\hfill$\Box$\end{trivlist}}

\begin{document}

	
\title{On isotropic points and contact invariants of curves}
\author{Amanda Dias Falqueto, Juan J. Nuño Ballesteros and Raúl Oset Sinha}

\maketitle
\begin{abstract}
Let $(C_{f},0) \subset (\mathbb{C}^{n},0)$ be a germ of an irreducible  holomorphic curve, given implicity by $f(x)=0$ if $n=2$ and by $f(x)=g(x)=0$ if $n=3$, and parametrised by $\gamma:(\mathbb{C},0) \rightarrow (\mathbb{C}^{n},0)$. We define the contact invariants $\operatorname{Iso}_{f}$ and $\operatorname{Iso}_{\gamma}$ for $f$ and $\gamma$, respectively, which count the number of isotropic points that appear after a generic deformation of the curve. We first consider the case where the curve is regular and establish a relation between those invariants and the number of vertices, inflection points and the zero-torsion points (for $n=3)$ that emerge from the deformation of the curve. Furthermore, we consider a germ of a singular curve with an isolated singularity and relate the number of isotropic points with the Milnor number of $f$.
\end{abstract}
	
	
	
\noindent \makebox [40mm]{\hrulefill}

\section{Introduction}\label{sec:intro}

We consider $\mathbb{C}^{n}$, with $n=2,3$, endowed with the complex bilinear form
\begin{align*}
    \langle .,. \rangle: \mathbb{C}^{n} &\rightarrow \mathbb{C}\\
    (v,w) &\mapsto \sum_{i=1}^{n}v_{i}w_{i},
\end{align*}
where $v=(v_{1}, \cdots, v_{n}), w = (w_{1}, \cdots , w_{n})$ (Pessers and Van der Veken \cite{Veken} call it the holomorphic inner product). We say that two vectors $v,w\in \mathbb{C}^{n}$ are orthogonal if $\langle v, w \rangle =0$ and a nonzero vector $v \in \mathbb{C}^{n}$ is isotropic if $\langle v, v \rangle =0$. Let $\gamma: D \rightarrow \mathbb{C}^{n}$ be a holomorphic curve, with $D$ an open subset of $\mathbb{C}$. We say that $t_{0} \in D$ is an isotropic point of  $\gamma$ if the tangent vector $\gamma'(t_{0})$ is an isotropic vector.
Isotropic points reveal new geometric information about the curve. For instance, Falqueto and Tari proved in \cite{FalquetoTari} that the evolute of a regular holomorphic curve admits a holomorphic extension at an isotropic point, where it is tangent to the original curve. Moreover, for an algebraic curve $X \subset \mathbb{C}^{2}$, Joiţa, Siersma and Tibar showed in \cite{Joitas} that the discriminant of a family of distance-squared functions, namely the Euclidean Discriminant (ED), can be decomposed into a non-linear component (the focal set of the curve) and other three linear components: the atypical discriminant, the singular discriminant and the iflex discriminant. 

In \cite{Struik}, Struik studied curves all of whose points are isotropic, called isotropic curves, and showed that the isotropic curves in $\mathbb{C}^{2}$ are the lines $x+iy = 0$ and $x-iy = 0$, namely the isotropic lines. In the complex 3-space, a plane is said to be isotropic if its normal vector is an isotropic vector. The normal to an isotropic plane at a point lies in this plane and it is tangent to the isotropic cone \textemdash meaning the cone generated by the isotropic lines \textemdash at that point.

For a germ of a holomorphic plane curve parametrised by $\gamma: (\mathbb{C},0) \rightarrow (\mathbb{C}^{2},0)$ and given by the zero set of a germ of a holomorphic function $f: (\mathbb{C}^{2},0) \rightarrow (\mathbb{C},0)$, we denote by $\operatorname{I}_{\gamma}$ and $\operatorname{V}_{\gamma}$ (resp. $\operatorname{I}_{f}$ and $\operatorname{V}_{f}$) the numbers of inflection points and vertices of the curve given by a parametrisation (resp. by an implicit equation). In \cite{FBM}, Bruce, Fernandes and Tari established a relation between the number of vertices and inflection points that appear after a generic deformation of a holomorphic singular curve in $\mathbb{C}^{2}$ {and proved that both are related to the type of the singularity of the curve. Additionally}, these geometric invariants depend on whether the curve is represented by a parametrisation or as the zero set of a map-germ, as the corresponding notions of generic deformation differ in these two settings.

By considering the isotropic points of a holomorphic curve with a non-isotropic limiting tangent line, we prove in this article that $\operatorname{I}_{\gamma}, \operatorname{I}_{f}, \operatorname{V}_{\gamma}$ and $\operatorname{V}_{f}$ are related to the number of isotropic points that emerge from a deformation of the curve and show that the latter depends only on the multiplicity of the curve. We also compute the number of isotropic points of a curve in terms of its irreducible branches, which allows us to restrict to the case where the curve is the zero set of an irreducible holomorphic function and to prove that the number of isotropic points of a curve is related to its Milnor number. 

 In \S \ref{S: regular}, we consider a regular holomorphic plane curve and compute the number of isotropic points that emerge under a deformation of the curve. We also establish a relation between this number and the number of inflection points and vertices of the curve.

In \S \ref{S: space}, we first consider an algebraic curve in $\mathbb{C}^{3}$ and compute the global number of its isotropic points, the number of points at which the osculating plane is isotropic, and the number of zero-torsion points. In the local setting, we compute how many of each type of these points arise after a deformation of the curve and we also establish a relation between them. We prove that, {both in the regular and in the singular case}, this relation depends on whether the limiting tangent direction of the curve at the singular point is isotropic or not. Indeed, if the limiting tangent direction of the curve at the singular point is an isotropic direction, the number of isotropic points depends on the contact between the curve and the isotropic cone. However, when the limiting tangent direction at the singular point is non-isotropic, this number is determined by the multiplicity of the curve.

\section{Preliminaries}

\begin{subsection}{Holomorphic map-germs}
Let $\mathcal{O}_{n,p}$ denote the  $\mathcal{O}_{n}$-module of germs of holomorphic maps $f:(\mathbb{C}^{n},0) \rightarrow \mathbb{C}^{p}$. If $p=1$, we denote it simply by $\mathcal{O}_{n}$. In this case, we define
\begin{align*}
    \mu(f) &= \dim_{\mathbb{C}} \frac{\mathcal{O}_{n}}{\langle f_{x_{1}}, \cdots, f_{x_{n}}\rangle}_{\mathcal{O}_{n}}
    \quad \mbox{ and }
\tau(f)= \dim_{\mathbb{C}}\frac{\mathcal{O}_{n}}{\langle f_{x_{1}}, \cdots, f_{x_{n}},f \rangle}_{\mathcal{O}_{n}}
 \end{align*}
as the Milnor Number of $f$ and the Tjurina number of $f$, respectively, where $f_{x_{i}}= \frac{\partial f}{\partial x_{i}}$ and $ \langle ,\rangle_{\mathcal{O}_{n}}$ denotes the ideal in $\mathcal{O}_{n}$ generated by the functions between the angle brackets. The former number counts the critical points of $f$, whereas the latter computes the number of parameters in a miniversal deformation of the zero set of $f$. 

The $\mathcal{K}$ group is the group of diffeomorphisms $H: (\mathbb{C}^{n} \times \mathbb{C}^{p},0) \rightarrow (\mathbb{C}^{n} \times \mathbb{C}^{p},0)$ such that $H(x,y)=(h(x), \theta(x,y))$, with $\theta(x,0)=0$ with $\theta \in \mathcal{O}_{{n+p},p}$ and $h \in \mathcal{O}_{n,n}$. If we set $h$ equal to the identity, then we have the $\mathcal{C}$ group. 
Two germs $f, g \in \mathcal{O}_{n,p}$ are $\mathcal{K}$-equivalent if there exist germs of diffeomorphisms $h: (\mathbb{C}^{n},0) \rightarrow (\mathbb{C}^{n},0)$ and $H:(\mathbb{C}^{n} \times \mathbb{C}^{p},0) \rightarrow (\mathbb{C}^{n} \times \mathbb{C}^{p},0)$, such that $H(x, y) = (h(x), \theta(x,y))$, with $\theta(x,0)= 0$, such that $H(x,f(x)) = (h(x), g(h(x)))$. If we consider $h$ as the germ of identity, the germs $f$ and $g$ are said to be $\mathcal{C}$-equivalent. 
A germ $f \in \mathcal{O}_{n}$ is k-$\mathcal{G}$-determined, where $\mathcal{G}= \mathcal{K}$ or $\mathcal{C}$, if for any germ $g \in \mathcal{O}_{n}$ such that $j^{k}g=j^{k}f$, we 
have $f$ $\mathcal{G}$-equivalent to $g$, where $j^{k}g$ is the $k$-jet of $g$ at the origin, that is, the Taylor polynomial of degree $k$ at the origin. We say that a germ $f$ is $\mathcal{G}$-finite if it is $k$-$\mathcal{G}$-finite for some integer $k$. If $\mathcal{G}= \mathcal{K}$, then $\mathcal{K}$-finiteness is equivalent to the germ having an isolated singularity at the origin.



Let $f = \sum a_{k}t^{k}$ be a germ of a holomorphic function. The order of $f$ is given by $
    \operatorname{ord}(f) := \min\{k; a_{k} \neq 0\}.
$
When $f = (f_{1}, \cdots, f_{k})$, $\operatorname{ord}(f) = \min_{\{i=1,\cdots, k\}}\operatorname{ord}(f_{i})$.
\end{subsection}

\subsection{Geometry of holomorphic curves}

Let $\gamma: D \rightarrow \mathbb{C}^{n}$ be a holomorphic curve, with $D \subset \mathbb{C}$ open and let $f=0$ be a hypersurface in $\mathbb{C}^{n}$, where $f: U \rightarrow \mathbb{C}$ is a holomorphic function and $U$ is an open subset of $\mathbb{C}^{n}$. We say that $\gamma$ and $f^{-1}(0)$ have order of contact $k$ at $t_{0} \in \mathbb{C}$ if $(f \circ \gamma)^{(i)}(t_{0})= 0$, for $i=0, \cdots, k-1$ and $(f \circ \gamma)^{(k)}(t_{0}) \neq 0$. The order of contact is at least $k$ if the last condition is omitted. 
A point where a plane curve has order of contact at least three with its tangent line is called an inflection point, whereas a vertex is a point where the contact of a plane curve and a circle is at least four.
If we consider holomorphic curves in $\mathbb{C}^{3}$, we say that $t_{0}$ is a zero-torsion point if the curve has contact at least four with its osculating plane at $t_{0}$. 

The contact of a holomorphic curve with model objects, such as lines, planes and circles, can be measured by the singularities of suitable map-germs and captures the local geometry of the curve, as in the real case (see \cite{FalquetoTari}, \cite{IRRT}). Indeed, the contact of a holomorphic plane (resp. space) curve with a line (resp. plane) orthogonal to a unit vector $v \in \mathbb{C}^{2}$ (resp. $v \in \mathbb{C}^{3}$) can be characterized by the singularity type of the height function, defined by
\begin{align*}
H_{v}: D &\rightarrow \mathbb{C}\\
t &\mapsto H_{v}(t) = \langle \gamma(t), v \rangle.
\end{align*}

In this setting, the origin is an inflection point of a holomorphic plane curve $\gamma$ if, for a vector $v$ orthogonal to $\gamma'(0)$, the function $H_{v}$ has an ${A}_{\geq 2}$-singularity at that point, that is, if ${det}(\gamma'(0), \gamma''(0))$ vanishes, where $\det$ denotes the determinant of the matrix. For $\gamma(t)= (x(t), y(t))$, this condition is equivalent to $i_{\gamma}(0)=0,$ where $i_{\gamma}:= x'y''- x''y'$. {Moreover, if $\gamma$ is the zero set of a holomorphic function $f: U \rightarrow \mathbb{C}$, where $U$ is an open subset of $\mathbb{C}^{2}$, by differentiating $f (x(t), y(t)) \equiv 0$ with respect to $t$ and expressing the derivatives of $x$ and $y$ in terms of the partial derivatives of $f$, $i_{\gamma}$ can be rewritten as} $i_{f}:=f_{xx}f_{y}^{2}-2f_{xy}f_{y}f_{x}+f_{yy}f_{x}^{2}$. Thus, the curve $f=0$ has an inflection point at the origin if $i_{f}$ vanishes at $t=0$.

In the complex space, the osculating plane of $\gamma$ at $t_{0}$ is the plane through $\gamma(t_{0})$ orthogonal to a vector $v \in \mathbb{C}^{3}$ such that $H_{v}$ has an $A_{\geq 2}$-singularity at $t_{0}$. Therefore, the orthogonal vector to the osculating plane is parallel to the vector $\gamma'(t) \times \gamma''(t)$, where $\times$ denotes the cross product. The origin is a zero-torsion point if $H_{v}$ has an ${A}_{\geq 3}$-singularity at that point, that is, if $det(\gamma'(0) \times \gamma''(0), \gamma'''(0))$ vanishes. 

We can describe the contact of a holomorphic curve $\gamma: D \rightarrow \mathbb{C}^{2}$ (resp. $D \rightarrow \mathbb{C}^{3})$, where $D$ is an open set of $\mathbb{C}$, with a circle centered at $p \in \mathbb{C}^{2}$ (resp. sphere centered at $ p \in \mathbb{C}^{3}$) by the singularities of the distance-squared function 
\begin{align*}
    D_{p}: D &\rightarrow \mathbb{C}\\
    t &\mapsto D_{p}(t) = \langle \gamma(t) -p, \gamma(t) -p\rangle.
\end{align*}
Thus, for $p = (a, b) \in \mathbb{C}^{2}$ and $\gamma(t) = (x(t), y(t))$, the origin is a vertex of the curve if $D_{p}$ has an ${A}_{\geq 3}$-singularity at that point. It follows that the origin is a solution of the system 
\begin{align*}
    x'(x-a)+y'(y-b) &=0\\
    x''(x-a)+y''(y-b)+x'^{2}+y'^{2}&=0\\
    x'''(x-a)+y'''(y-b)+3x''x'+3y''y'&=0.
\end{align*}
Equivalently, using the matrix representation of the system above, we have that the origin is a vertex of the curve if the determinant of the matrix 
$$\begin{pmatrix}
xx'+yy' & x' & y'\\
x'^{2}+y'^{2}+xx''+yy'' & x'' & y''\\
3x'x''+3y'y''+xx'''+yy''' & x''' & y'''
\end{pmatrix}$$
is zero at that point, that is, if $v_{\gamma}(0)=0$, where 

$$v_{\gamma}:= (x'^{2}+y'^{2})(x'y'''-x'''y')+3(x'x''+y'y'')(x''y'-x'y'').$$

The tangency of two germs of irreducible plane curves given by the zero sets of $f, g \in \mathcal{O}_{2}$ 
can be measured by the intersection multiplicity, which is defined by
\begin{align*}
    m(f,g) &= \dim_{\mathbb{C}} \frac{\mathcal{O}_{2}}{\langle f,g \rangle}_{\mathcal{O}_{2}}\\
    &= \operatorname{ord}(f \circ \gamma),
\end{align*}
where $\gamma$ is the parametrisation of the curve defined by $\{g=0\}$ (see the properties of the intersection multiplicity in \cite{Gruel}).

\section{Isotropic points on plane curves}\label{S: regular}

Let $\gamma: D \rightarrow \mathbb{C}^{2}, \gamma(t)= (x(t),y(t))$, with $D \subset \mathbb{C}$ open, be a parametrisation of a holomorphic curve $C_{f}$ given by $\{f=0\}$, where $f: U \rightarrow \mathbb{C}$ is a holomorphic function and $U$ is an open subset of $\mathbb{C}^{2}$. The isotropic points of the curve $\gamma$ are the zeros of
\begin{align*}
    \operatorname{iso}_{\gamma}:= x'^{2}+y'^{2},
\end{align*} 
and the number of isotropic points of a deformation of the curve $\gamma$, counted with multiplicity, is given by
\begin{align*}
    \operatorname{Iso}_{\gamma}= \dim_{\mathbb{C}}\frac{\mathcal{O}_{1}}{\langle \operatorname{iso}_{\gamma} \rangle}_{\mathcal{O}_{1}}.
\end{align*}
We say that the origin is an isotropic point of $C_{f}$ if the origin is a solution of
\begin{align*}
    f(x,y)&=0 \\
    \operatorname{iso}_{f}=f_{x}(x,y)^{2}+f_{y}(x,y)^{2}&=0, 
\end{align*}
and the number of isotropic points, counted with multiplicity, on a deformation of the curve $C_{f}$ is
\begin{align*}
    \operatorname{Iso}_{f}= \dim_{\mathbb{C}}\frac{\mathcal{O}_{2}}{\langle f, \operatorname{iso}_{f}\rangle}_{\mathcal{O}_{2}}.
\end{align*}

\begin{rem}
    {It follows from the definition that a singular point is an isotropic point.}
\end{rem}

\begin{defn}
    Let $U \subset \mathbb{C}^{2}$ be open. A holomorphic map  $\phi: U \rightarrow \mathbb{C}^{2}$ is conformal if for all $p \in U$ and for all $v_{1}, v_{2} \in \mathbb{C}^{2}$ we have
    \begin{align*}
        \langle d\phi_{p}(v_{1}), d\phi_{p}(v_{2}) \rangle = \alpha(p) \langle v_{1}, v_{2} \rangle,
    \end{align*}
    where $\alpha: U \rightarrow \mathbb{C}$ is a holomorphic function that does not vanish in $U$.
\end{defn}
\begin{prop} Let $\gamma: (\mathbb{C},0) \rightarrow (\mathbb{C}^{2},0)$ be a holomorphic curve defined implicitly by $\{f=0\}$, where $f:(\mathbb{C}^{2},0) \rightarrow (\mathbb{C},0)$ is a germ of a holomorphic function. Then, 
\begin{enumerate}
    \item The integer $\operatorname{Iso}_{f}$ is a $\mathcal{C}$-invariant of $f$, that is, invariant under $\mathcal{C}$-equivalence.
    \item Assume that $f$ has an isolated isotropic point at the origin and let $F: U \times V,0  \rightarrow \mathbb{C}, F(u,x,y) = f_{u}(x,y)$, where $U$ and $V$ are open neighbourhoods of the origin in $\mathbb{C}^{s}$ and in $\mathbb{C}^{2}$, respectively, and $F(0,x,y) = f_{0}(x,y)=f(x, y)$. Then,
\begin{align*}
\operatorname{Iso}_{f} = \sum_{p \in V(f_{u}, \operatorname{iso}_{f_{u}})\cap V} \operatorname{Iso}_{f_{u}}(p).
   \end{align*}
\item $\operatorname{Iso}_{f} \geq \tau(f)$, where $\tau(f)$ is the Tjurina number of $f$.
\item Let $\gamma:(\mathbb{C},0) \rightarrow (\mathbb{C}^{2},0), \gamma(t) = (x(t), y(t))$ and let $\phi: (\mathbb{C}^{2},0) \rightarrow (\mathbb{C}^{2},0)$ be a conformal mapping. Then, the curves $\gamma$ and $\beta$ have the same number of isotropic points, where $\beta = \phi \circ \gamma$.
\end{enumerate}
\end{prop}
\begin{proof}
\begin{enumerate}
\item Let $f, g \in \mathcal{O}_{2}$ such that $f$ is $\mathcal{C}$-equivalent to $g$. In particular,  there exists $h \in \mathcal{O}_{2}$ with $h(0,0) \neq 0$ such that $g= hf$. Differentiating $g$ with respect to $x$ and $y$ we get
\begin{align*}
    g_{x}&=h_{x}f+f_{x}h \\
     g_{y}&=h_{y}f+f_{y}h.
\end{align*}

Hence,
\begin{align*}
    g_{x}^{2}+ g_{y}^{2} &= f(h_{x}^{2}f+ 2f_{x}h_{x}h+ 2f_{y}hh_{y}+fh_{y}^{2})+h^{2}(f_{x}^{2}+ 
      f_{y}^{2}).
\end{align*}

Therefore, $ \langle g, g_{x}^{2}+ g_{y}^{2}\rangle_{\mathcal{O}_{2}}  \subset \langle f, f_{x}^{2}+f_{y}^{2}\rangle_{\mathcal{O}_{2}} $ and, thus, $\operatorname{Iso}_{g} \geq \operatorname{Iso}_{f}$. 
Using $f= h^{-1}g$, we obtain the reverse inequality. 

\item Consider the ideal $\langle F, {F}_{x}^{2}+{F}_{y}^{2} \rangle$ in $\mathcal{O}_{2+s}$. Since the origin is an isolated isotropic point of $F$, we have $\dim V(F, F_{x}^{2}+F_{y}^{2}) \geq s$. Considering the projection
\begin{align*}
    \pi: V(F, F_{x}^{2}+F_{y}^{2}) \rightarrow (\mathbb{C}^{s},0)
\end{align*}
we get that $\pi^{-1}(0) =  V(f, f_{x}^{2}+f_{y}^{2}) = \{0\}$ is finite hence $\pi$ is a finite morphism. It follows from Remmert's Finite Mapping Theorem (see \cite{Loja}, Chapter V) that $\pi(V(F, F_{x}^{2}+F_{y}^{2}))$ is an analytic subset of $(\mathbb{C}^{s}, 0)$ and has the same dimension of $V(F,F_{x}^{2}+F_{y}^{2})$. So $\dim V(F, F_{x}^{2}+F_{y}^{2}) \leq s$. Thus, $\dim V(F, F_{x}^{2}+F_{y}^{2})=s$. 
The ring $\frac{\mathcal{O}_{2+s}}{\langle F, F_{x}^{2}+F_{y}^{2}\rangle}_{\mathcal{O}_{2+s}}$ is a complete intersection and hence, Cohen-Macaulay. Therefore, the principle of conservation of multiplicity applies.
\item Since
\begin{align*}
    f_{x}^{2}+f_{y}^2={(f_{x}+if_{y})(f_{x}-if_{y})},
\end{align*}
we have the inclusion
\begin{align*}
    \langle f_{x}^{2}+f_{y}^{2} \rangle_{\mathcal{O}_{2}} \subset \langle f_{x}-if_{y},f_{x}+if_{y} \rangle_{\mathcal{O}_{2}} = \langle f_{x}, f_{y}\rangle_{\mathcal{O}_{2}}.
\end{align*}
 We then obtain
\begin{align*}
\dim \frac{\mathcal{O}_{2}}{\langle f, f_{x}^{2}+f_{y}^{2} \rangle}_{\mathcal{O}_{2}} \geq \dim \frac{\mathcal O_{2}}{\langle f, f_{x}, f_{y}\rangle}_{\mathcal{O}_{2}}= \dim \frac{\mathcal O_{2}}{J_{f} + \langle f \rangle}_{\mathcal{O}_{2}} = \tau(f).
\end{align*}

\item Let $\phi: (\mathbb{C}^{2},0) \rightarrow (\mathbb{C}^{2},0)$ be a conformal mapping, that is, $\langle d\phi_{p} v_{1}, d\phi_{p}v_{2} \rangle = \alpha(p)\langle v_{1}, v_{2}\rangle$, for some germ of a holomorphic function  $\alpha: (\mathbb{C}^{2},0) \rightarrow \mathbb{C}$ and consider $\beta:= \phi \circ \gamma$. The result follows from the following equality
\begin{align*}
\langle \beta '(t), \beta'(t) \rangle = \langle d\phi_{\gamma(t)}(\gamma'(t)), d\phi_{\gamma(t)}(\gamma'(t))\rangle = \alpha(\gamma(t))\langle\gamma'(t), \gamma'(t)\rangle . 
\end{align*}
\end{enumerate}
\end{proof}

We can also define an isotropic point of a plane curve by observing its contact with the isotropic lines. In this context, an isotropic point of a regular curve is a point where the curve has order of contact at least two with one of the two isotropic lines through it. We study the isotropic points considering the isotropic line $x+iy=0$. The same hold if we consider the other isotropic line.  
We first compute the number of isotropic points for generic germs of nonsingular curves in $\mathbb{C}^{2}$, given as the zero set of a germ of a holomorphic function. 

\subsection{Regular curves}

\begin{theo}
Let $C_{f}$ be a germ of an irreducible plane curve, with $f = x+iy+g_{k} + O(k+1)$, where $g_{k}$ is a homogeneous polynomial of degree $k$ and $k \geq 2$. Then, $C_{f}$ has $k-1$ isotropic points, provided the line $x+iy=0$ is not contained in the tangent cone of ${g_{{k}}}_{x}+ i{g_{{k}}}_{y}=0$. 
\end{theo}

\begin{proof}
Note that
\begin{equation}
\begin{aligned}
   (f_{x}+if_{y})&= (1+{g_{k}}_{x} + \mbox{h.o.t.})+i(i+{g_{k}}_{y} + \mbox{h.o.t.})\\
   &={g_{k}}_{x}+i{g_{k}}_{y} + \mbox{h.o.t.}
\end{aligned}
\end{equation}
Factoring the leading homogeneous part of $f_{x}+if_{y}$ we obtain
\begin{align*}
    g_{{k}_{x}}+ig_{{k}_{y}} = \prod_{j = 1}^{k-1}(\alpha_{j}x+\beta_{j}y).
\end{align*}

Then,  
\begin{align*}
   \operatorname{Iso}_{f}&= \dim_{\mathbb{C}} \frac{\mathcal O_{2}}{\langle f, f_{x}+if_{y} \rangle}_{\mathcal{O}_{2}}\\
    &= \dim_{\mathbb{C}} \frac{\mathcal O_{2}}{\langle x+iy + O(k), \prod_{j=1}^{k-1}(\alpha_{j}x+\beta_{j}y)+ O(k) \rangle}_{\mathcal{O}_{2}}\\
    &=  \dim_{\mathbb{C}}  \frac{\mathcal O_{2}}{\langle x+iy, \prod_{j=1}^{k-1}\alpha_{j}x+\beta_{j}y\rangle}_{\mathcal{O}_{2}}\\
    &= k-1,
    \end{align*}
where the third equality follows from {\cite{teojuanjo}, Section 3.3, Corollary 3}, provided that the determinants of the matrices

\[
B_{j} = \begin{bmatrix}
  1   & i \\
\alpha_{j} & \beta_{j}
\end{bmatrix}
\]
are nonzero, for all
$j=1, \cdots, k-1$. This condition means that the tangent cone of ${g_{k}}_{x}+i{g_{k}}_{y} = 0$ is not parallel to the vetor $(1,i)$. 
\end{proof}

In the following theorem we assume that the curve has an isotropic point and an inflection at the origin, and establish a relation between isotropic points and inflections for germs of regular holomorphic plane curves.

\begin{theo}\label{inflection} 
Let $\gamma:(\mathbb{C},0) \rightarrow (\mathbb{C}^{2},0), \gamma(t) = (x(t), y(t))$ be a germ of a regular holomorphic curve with an isotropic point and an inflection at the origin. Suppose that $\operatorname{Iso}_{\gamma}$ is finite. Then, $\operatorname{I}_{\gamma} = \operatorname{Iso}_{\gamma}-1$.
\end{theo}
\begin{proof}
Since $\gamma$ is regular, it can be locally parametrised as $\gamma(t) = (t, a_{1}t+O(2))$. It follows from the fact that $\operatorname{Iso}_{\gamma}$ is finite that there exists $j>1$ such that $a_{j} \neq 0$. We then write $\gamma(t)  =  (t, a_{1}t+a_{k}t^{k} + \mbox{h.o.t.})$, with $k = \min \{i; a_{i} \neq 0\}$. 

Since the origin is an isotropic point, we conclude that $a_{1} = \pm i$. Thus, 
     \begin{align}\label{iso}
         \operatorname{Iso}_{\gamma} = \operatorname{ord}(x'(t)^{2}+y'(t)^2) = \operatorname{ord}(\pm 2ia_{k}kt^{k-1} + \mbox{h.o.t.}) = k-1.
     \end{align}

Moreover, since the origin is an inflection point, we have $k \geq 3$ and 
\begin{align}\label{inf}
\mbox{I}_{\gamma} = \operatorname{ord}(det(\gamma',\gamma'')) = \operatorname{ord}(k(k-1)a_{k}t^{k-2} + \mbox{h.o.t.}) =k-2.
\end{align}

The result follows from \eqref{iso} and \eqref{inf}.
\end{proof}

The number of isotropic points of a germ of a regular holomorphic curve is also related to the number of vertices that emerge after a deformation of the curve, as shown in the following theorem.

\begin{theo}\label{vertice} 
Let $\gamma:(\mathbb{C},0)\rightarrow (\mathbb{C}^{2},0), \gamma(t) = (x(t), y(t))$ be a germ of a regular holomorphic curve with an isotropic point and a vertex at the origin. {Suppose that $\operatorname{Iso}_{\gamma}$ is finite}. Then, 
    \begin{align*}
\operatorname{V}_{\gamma} = 2\operatorname{I}_{\gamma} = 2(\operatorname{Iso}_{\gamma} - 1).
    \end{align*}
\end{theo}
\begin{proof}
Let $\gamma(t) = (t, it+a_{k}t^{k} + \mbox{h.o.t.})$ be a parametrisation of the curve, where $k= \min\{i; a_{i} \neq 0\}$.
 Since the origin is a vertex of $\gamma$, we have $k \geq 3$ and the distance-squared function has an ${A}_{\geq 3}$-singularity at that point. Hence, we have $v_{\gamma}(0) = 0$, where
\begin{align*}
    v_{\gamma} =(x'y'''-x'''y')(x'^{2}+y'^{2})-3(x'x''+y'y'')(x'y''-x''y').
\end{align*}
Substituting the derivatives of $\gamma$ in the expression above yields
\begin{align*}
v_{\gamma} =-ik^{2}(k^{2}-1)a_{k}^{2}t^{2k-4} + \mbox{h.o.t.}
\end{align*} 
Hence,
\begin{align*}
   \mbox{V}_{\gamma}= \dim_{\mathbb{C}}  \frac{\mathcal
{O}_{1}}{\langle v_{\gamma}\rangle}_{\mathcal{O}_{1}}= \dim_{\mathbb{C}}  \frac{\mathcal{O}_{1}}{\langle t^{2k-4} \rangle}_{\mathcal{O}_{1}}= 2k-4.
\end{align*}
It follows by \eqref{iso} and \eqref{inf} that
\begin{align*}
    \operatorname{V}_{\gamma} = 2 \operatorname{I}_{\gamma} = 2 (\operatorname{Iso}_{\gamma} -1).
\end{align*}

\end{proof}
\begin{rem}

It follows from the expressions for $v_{\gamma}$ and $i_{\gamma}$ that, for a regular holomorphic curve with an isotropic point at the origin, the origin is a vertex if and only if it is an inflection point. Then, the hypotheses of Theorem \ref{vertice} and of Theorem \ref{inflection} coincide.
\end{rem}

\begin{rem}
We observe that, generically, a point of a holomorphic curve cannot be simultaneously an isotropic point and an inflection point. Hence, the numbers of inflections and isotropic points arising under a deformation of the curve are well-defined. By the Basic Lema of Transversality \cite{Ballesteros}, a generic holomorphic curve does not have a point that is simultaneously an inflection and an isotropic point.
\end{rem}

In fact, let $\gamma: D \rightarrow \mathbb{C}^{2}, \gamma(t) = (t, it+ a_{1}t+ \cdots+ a_{k}t^{k}$) be a holomorphic curve, where $D$ is an open subset of $\mathbb{C}$, and define the set $\Sigma := \{ (x(t), y(t)) \in \mathbb{C}^{2}; x'_{a}(t)^{2}+y'_{a}(t)^{2} = x'_{a}(t)y''_{a}(t)-y'_{a}(t)x''_{a}(t)=0, \mbox{ for } a \in \mathbb{C}^{k}\}$. Consider the map
\begin{align*}
   F: D \times \mathbb{C}^{k} &\rightarrow \mathbb{C}^{4}\\
   (t, a) &\mapsto (x_{a}'(t), y_{a}'(t), x_{a}''(t), y_{a}''(t)).
\end{align*}
Since $x'_{a}\equiv 1$ and $x''_{a} \equiv 0$ for all $a \in \mathbb{C}^{k}$, we can rewrite the map $F$ as follows
\begin{align*}
     F: D \times \mathbb{C}^{k} &\rightarrow \mathbb{C}^{2}\\
    (t, a) &\mapsto (y_{a}'(t), y_{a}''(t)).
\end{align*}

We show that $F \pitchfork \Sigma$; then, by the Basic Lema of Transversality, $F_{a} \pitchfork \Sigma$, for almost every $a \in \mathbb{C}^{k}$.

Since $F(t, a) = ((a_{1}+i) + 2a_{2}t + 3a_{3}t^{2} + \mbox{h.o.t.}, 2a_{2} + \mbox{h.o.t.})$ and the minor of the Jacobian matrix of $F$ given by
\begin{align*}
  \begin{pmatrix}
 1 & 2t\\
    0 & 2
\end{pmatrix}
\end{align*}
has nonzero determinant, $F$ is a submersion and the result follows.

\subsection{The number of isotropic points of singular plane curves}
We consider a germ of a singular plane curve parametrised by $\gamma:(\mathbb{C},0) \rightarrow (\mathbb{C}^2,0)$, {which is an injective immersion away from the origin} and given implicitly by $f=0$, where $f:(\mathbb{C}^{2},0) \rightarrow (\mathbb{C},0)$ is a holomorphic function {with an isolated singularity at the origin}. We then compute the number of isotropic points, vertices and inflections that appear after a deformation of each curve, relating them to the Milnor number of the curve. 

\begin{theo}
Let $f=gh$.
    \begin{enumerate}
    \item (\cite{FBM}) If $\mbox{I}_{g}$ and $\mbox{I}_{h}$ are finite, then $\mbox{I}_{f}$ is finite and 
    \begin{align*}
        \mbox{I}_{f} = \mbox{I}_{g} + \mbox{I}_{h} + 6m(g,h).
    \end{align*}
    \item (\cite{FBM}) If $\mbox{V}_{g}$ and $\mbox{V}_{h}$ are finite, then $\mbox{V}_{f}$ is finite and 
    \begin{align*}
        \mbox{V}_{f} = \mbox{V}_{g} + \mbox{V}_{h} + 12 m(g,h).
    \end{align*}
    \item If $\operatorname{Iso}_{g}$ and $\operatorname{Iso}_{h}$ are finite, we have $\operatorname{Iso}_{f}$  finite and 
    \begin{align*}
        \operatorname{Iso}_{f} = \operatorname{Iso}_{g}+\operatorname{Iso}_{h} + 4m(g,h).
    \end{align*}
    \end{enumerate}
\end{theo}
\begin{proof}
    $(3)$ We observe that $\operatorname{iso}_{f} = g_{x}^{2}h^{2} + 2h_{x}hgg_{x}+ h_{x}^2g^{2}+g_{y}^{2}h^{2}+2g_{y}h_{y}gh+h_{y}^{2}g^{2}$. By the properties of the intersection number of $g$ and $h$ (see \cite{multiplicity}), we obtain
    \begin{align*}
        \operatorname{Iso}_{f} = & m(gh, g_{x}^{2}h^{2} + 2h_{x}hgg_{x}+ h_{x}^2g^{2}+g_{y}^{2}h^{2}+2g_{y}h_{y}gh+h_{y}^{2}g^{2})\\
        &= m(g,h^{2}\operatorname{iso}_{g}) + m(h, g^{2}\operatorname{iso}_{h})\\
        &= \operatorname{Iso}_{g} + \operatorname{Iso}_{h} + 4m(g,h).
    \end{align*}
\end{proof}

The preceding theorem yields by induction:
\begin{theo} \label{teo310}
    If $f = f_{1} \cdots f_{n}$, then
 $\operatorname{Iso}_{f} = \sum_{i=1}^{n} \operatorname{Iso}_{f_{i}} + 4 \sum_{i<j}m(f_{i},f_{j}).$
\end{theo}

In \cite{FBM}, it was shown that 
the numbers of vertices and inflection points that appear after a deformation of a holomorphic curve $f=0$ are completely determined by the irreducible components of $f$. Thus, in order to determine the number of inflections, vertices and isotropic points concentrated at the singularity, it is enough to consider $f$ irreducible.

Therefore, we consider a germ of an irreducible curve $C_{f}$ given by $\{f=0\}$, with multiplicity $m$ and its parametrisation given by $\gamma:(\mathbb{C},0) \rightarrow (\mathbb{C}^{2},0), \gamma(t) = (t^{m}, y(t))$, with $y \in \mathcal{O}_{1}$ and $\operatorname{ord}(y)>m$. We establish a relation between the numbers of isotropic points $\operatorname{Iso}_{\gamma}$ and $\operatorname{Iso}_{f}$ and the Milnor number $\mu_{f}$ in an analogous way to \cite{FBM}, where the authors showed that $\operatorname{I}_{f} = \operatorname{I}_{\gamma} + 3 \mu(f)$ and 
$\operatorname{V}_{f}= \operatorname{V}_{\gamma} + 6 \mu(f)$.

\begin{theo}\label{usarei} Let $C_{f}$ be a germ of an irreducible holomorphic curve and let 
 $\gamma: (\mathbb{C},0) \rightarrow (\mathbb{C}^{2},0)$ be its parametrisation. Then,
 $\operatorname{ Iso}_{f} = \operatorname{ Iso}_{\gamma} + 2 \mu(f).$
 \end{theo}
\begin{proof}
   Differentiating $f(x(t), y(t))=0$ with respect to $t$ yields
    \begin{align*}
        f_{x}x'+ f_{y}y' = 0.
    \end{align*}
Therefore, there exists a germ of a holomorphic function $h: (\mathbb{C},0) \rightarrow \mathbb{C}$, with $h(0) \neq 0$ such that 
\begin{align}
f_{x}(x(t),y(t)) &= h(t)t^{k}y'(t) \label{ordemx}\\
f_{y}(x(t), y(t)) &= -h(t)t^{k}x'(t).\label{ordemy}
\end{align}
 
Since $m$ is the multiplicity of the curve, we have 
\begin{align*}
    f(0,y) = \frac{1}{m!}\frac{\partial^{m} f}{\partial y^{m}}(0,0) y^{m} + \mbox{ h.o.t.}
\end{align*}
As $f(0,y)$ does not vanish for any $y\neq 0$ in a neighbourhood of the origin, it follows by Teissier's Lemma (see \cite{Teissier}) that \begin{align}\label{Milnor}
    m(f(x,y), f_{y}(x,y)) &= \mu(f) + m(f(x,y),x) -1\nonumber\\
    &= \mu(f) + \operatorname{ord} (f(0,y))-1 \nonumber \\
    &= \mu(f)+ m-1.
\end{align}
From \eqref{ordemy}, we obtain
\begin{align}\label{ord y}
    \operatorname{ord}(f_{y}(x(t),y(t)))= k+\operatorname{ord}(x'(t)) = k+m-1 
\end{align}

It follows from \eqref{Milnor} and \eqref{ord y} that
\begin{align}\label{milnornumber}
    \mu(f) = k.
\end{align}

Using \eqref{ordemx} and \eqref{ordemy}, we have:

\begin{align*}
    (f_{x}^{2}+f_{y}^2) (x(t), y(t)) &= h(t)^{2}t^{2k}y'(t)^{2}+ h(t)^{2}t^{2k}x'(t)^{2}\\
    &= h(t)^{2}t^{2k}(x'(t)^{2}+y'(t)^{2})\\
    &= h(t)^{2}t^{2k}\operatorname{iso}_{\gamma}.
\end{align*}

Hence,
\begin{align}
    \operatorname{ Iso}_{f} = m(f,\operatorname{iso}_{f}) &= \operatorname{ord}(h(t)^{2}t^{2k}\operatorname{iso}_{\gamma}) \nonumber \\
    &= \operatorname{ord}(\operatorname{iso}_{\gamma})+2k. \label{igualdadeiso}
\end{align}

It follows from \eqref{milnornumber} and \eqref{igualdadeiso} that
\begin{align*}
    \operatorname{Iso}_{f} = 2\mu(f) + \operatorname{Iso}_{\gamma}.
\end{align*}
\end{proof}

For $\gamma:(\mathbb{C},0) \rightarrow (\mathbb{C}^{2},0), \gamma(t) = (t^{m}, at^{n}+O(n+1)), n>m$, there exists a unique circle with maximal order of contact with the curve. Denoting this order of contact by $\lambda(\gamma)$, the authors in \cite{FBM} showed that  $\mbox{V}_{\gamma}= \mbox{I}_{\gamma} + \lambda(\gamma) -3$. Moreover, the number of inflection points that arise under a deformation of the curve is given by $\mbox{I}_{\gamma}= m+n-3$. We shall need these results and the following theorem.

\begin{theo}\label{tudo}
Let $\gamma: (\mathbb{C},0) \rightarrow (\mathbb{C}^{2},0)$ be a germ of a holomorphic curve with a non-isotropic limiting tangent direction of the curve at the origin. So, we may consider the parametrisation by $\gamma(t) = (t^{m}, at^{n}+O(n+1))$, with $n>m$ and $a \neq 0$. Then,
\begin{enumerate}\label{teoremavertice}
    \item \label{isodegamma} $\operatorname{Iso}_{\gamma}= 2(m-1)$; 
    \item $\operatorname{V}_{\gamma} =  \frac{1}{2}\operatorname{Iso}_{\gamma}+\lambda(\gamma)+n-5$.
    \end{enumerate}
\end{theo}






\begin{proof}
\begin{enumerate}
   \item  Note that the curve has contact of order at least two with an isotropic line if and only if $x'(t) +iy'(t) =0$ or $x'(t) -iy'(t) =0$. 
Since the order of  $x'+ iy'$ and $x'-iy'$ are both equal to $m-1$, we have
$\operatorname{Iso}_{\gamma} = 2(m-1)$.

\item The proof follows combining $\operatorname{I}_{\gamma}= m+n-3$, $\operatorname{V}_{\gamma} = \operatorname{I}_{\gamma}+\lambda(\gamma)-3$ and item 1.
\end{enumerate}
\end{proof}

\begin{rem}
{For a germ of a holomorphic curve parametrised by }$\gamma: (\mathbb{C},0) \rightarrow (\mathbb{C}^{2},0), \gamma(t) = (t^{m}, it^{m}+at^{n}+O(n+1))$, that is, whose limiting tangent direction at the origin is isotropic, we have $\operatorname{Iso}_{\gamma} = \operatorname{ord}(x'(t)^2+y'(t)^2) = \operatorname{ord } (\pm 2aimnt^{m+n-2} + \mbox{h.o.t.}) = m+n-2.$   

We observe that, although the number of isotropic points $\operatorname{Iso}_{\gamma}$ depends on whether the limiting tangent direction is isotropic, the equality $\operatorname{ Iso }_{f} = \operatorname{Iso}_{\gamma}+ 2\mu(f)$ holds independently from this condition, since the proof of Theorem \ref{usarei} only requires the curve to have multiplicity $m$.
\end{rem}

If $C_{f}$ is a germ of a singular irreducible curve {with a non-isotropic tangent line at the origin} and $\lambda(f)$ is the maximal order of contact between $C_{f}$ and circles at the origin, we have 
    $\operatorname{V}_{f}=\operatorname{I}_{f}+3\mu(f)+ \lambda(f).$ (see \cite{FBM}).

The next theorem shows the relation between the number of isotropic points, the geometric invariants of a singular plane curve and its Milnor number.

\begin{theo}
Let $C_{f}$ be a germ of a singular  
irreducible curve whose order of contact with its non-isotropic tangent line at the origin is $n$. Then,
\begin{align*}
\operatorname{Iso}_{f} = \operatorname{V}_{f}-\lambda(f)- 7\mu(f)+\operatorname{I}_{f}-2n+4.
\end{align*}
\end{theo}
\begin{proof}
 It follows by Theorem \ref{usarei} and by $\operatorname{I}_{f} = \operatorname{I}_{\gamma}+3\mu(f)$ that 
 \begin{align*}
    \operatorname{Iso}_{f}-\operatorname{I}_{f}&= \operatorname{Iso}_{\gamma}-\operatorname{I}_{\gamma} -\mu(f)\\
     &= \operatorname{I}_{\gamma}- \mu(f)-2n+4.
    \end{align*}
Since $\operatorname{I}_{f} = \operatorname{I}_{\gamma} + 3 \mu(f)$ and $\operatorname{V}_{f} = \operatorname{I}_{f}+3 \mu(f)+\lambda(f)$, we have
\begin{align*}
\operatorname{Iso}_{f} &= 2\operatorname{I}_{f}-4\mu(f)-2n+4\\
    &= \operatorname{I}_{f} +3\mu(f)-7\mu(f)+\operatorname{I}_{f}-2n+4\\
   &= \operatorname{V}_{f}-\lambda(f)- 7\mu(f)+\operatorname{I}_{f}-2n+4.
\end{align*}
\end{proof}

\section{Space curves}\label{S: space}

\subsection{Geometric invariants of holomorphic space curves}

We first consider algebraic curves in $\mathbb{C}^{3}$ and compute the global number of isotropic points, zero-torsion points and points where the osculating plane is isotropic. The number of isotropic points is defined in an analogous way to that of the plane curves.

\begin{prop}Let $f, g \in \mathbb{C}[x,y,z]$ be two polynomials of degrees $d_{1}$ and $d_{2}$, respectively, and let  $F(x:y:z:w) = f_{d_{1}}(x,y,z)+wf_{d_{1-1}}(x,y,z) + \cdots + w^{d_{1}}f_{0}(x,y,z)$ and $G(x:y:z:w) = g_{d_{2}}(x,y,z) + wg_{d_{2-1}}(x,y,z)+ \cdots + w^{d_{2}}g_{0}(x,y,z)$ be the homogenizations of $f$ and $g$. Consider the algebraic curve $C_{\Phi}= \Phi^{-1}(0), \Phi =(f,g)$, given by the intersection of the surfaces $f=0$ and $g=0$. Assume that the projective closure of the curve $\Phi$ is regular at every point at infinity and intersects the plane $w=0$ transversely. If the system
\begin{align*}
    {f}_{{d}_{1}}(\alpha) &= 0\\
     {g}_{{d}_{2}}(\alpha) &=0\\
     \langle \alpha , \alpha \rangle &=0
\end{align*}
has no solution in the plane at infinity $\{w=0\}$ and the number of isotropic points of $\Phi$ is finite, the curve $C_{\Phi}$ has $\operatorname{Iso}_{\Phi}=2d_{1}d_{2}(d_{1}+d_{2}-2)$ isotropic points in $\mathbb{C}^{3}$.
\end{prop}
\begin{proof}
 Let $\gamma: D \rightarrow \mathbb{C}^{3}$, $\gamma(t) = (x(t), y(t), z(t))$ be a parametrisation of the curve $\Phi$, with $D$ an open subset of $\mathbb{C}$. Differentiating $\Phi(\gamma(t)) = 0$ with respect to t yields
\begin{align*}
    d\Phi(\gamma(t)) \cdot \gamma' = 0,
\end{align*}
which implies that $\gamma'(t) = \nabla f(\gamma(t)) \times \nabla g(\gamma(t))$. 

The isotropic points of the curve $f=g=0$ are the solutions of the system
\begin{equation}\label{isotropic}
\begin{cases}
    f &= 0\\
    g &= 0\\
    (f_{y}g_{z}-f_{z}g_{y})^{2}+(f_{z}g_{x}-g_{z}f_{x})^{2}+(f_{x}g_{y}-f_{y}g_{x})^{2} &=0
\end{cases}
\end{equation}

Homogenising the polynomials $f$ and $g$, we obtain the homogenisation of the third polynomial in the above system. Then, by Bézout's theorem, the system \eqref{isotropic} has $2d_{1}d_{2}(d_{1}+d_{2}-2)$ solutions in $\mathbb{CP}^{3}$, counted with multiplicity. 

Let $\alpha = (\alpha_{1}:\alpha_{2}:\alpha_{3}:0) \in \mathbb{CP}^{3}$ be a solution at infinity of the above system. From Euler's identity, we have
\begin{align}
    \alpha_{1}f_{{d_{1}}_{x}}(\alpha) + \alpha_{2}f_{{d_{1}}_{y}}(\alpha) + \alpha_{3}f_{{d_{1}}_{z}}(\alpha) &=0\label{fdx}\\
    \alpha_{1}g_{{d_{2}}_{x}}(\alpha) + \alpha_{2}g_{{d_{2}}_{y}}(\alpha) + \alpha_{3}g_{{d_{2}}_{z}}(\alpha) \label{gdx} &=0
\end{align}


Assume $\alpha_{1} \neq 0$. Substituting the expressions of ${f_{d_{1}}}_{x}(\alpha)$ and ${g_{d_{2}}}_{x}(\alpha)$ into the third equation of \eqref{isotropic} yields
\begin{align}
    (1+(\frac{\alpha_{2}}{\alpha_{1}})^{2}+ (\frac{\alpha_{3}}{\alpha_{1}})^{2})( {f_{d_{1}}}_{y}{g_{d_{2}}}_{z} - {f_{d_{1}}}_{z}{g_{d_{2}}}_{y})^2=0,
\end{align}

which implies that $\alpha$ is an isotropic point at infinity  or $({f_{d_{1}}}_{y}{g_{d_{2}}}_{z} - {f_{d_{1}}}_{z}{g_{d_{2}}}_{y})^2=0$. The last condition is equivalent to the surfaces $F=0$ and $G=0$ intersecting tangentially at $w=0$, or the projective closure of $\Phi$ having a singularity at $\alpha$.

As both are excluded by hypothesis, it follows that the projectivized system has no roots at infinity. Therefore, the curve has $2d_{1}d_{2}(d_{1}+d_{2}-2)$ isotropic points in $\mathbb{C}^{3}$.
\end{proof}

\begin{rem}
The global number of isotropic points can also be calculated using the  General Plucker Formula. Indeed, since the isotropic vector of the curve is represented as a direction in the conic $\{x^{2}+y^{2}+z^{2}=0, w=0\}$, the number of isotropic points of the curve in $\mathbb{P}^{3}$ coincides with the number of points of the curve where the tangent line intersects that conic, which is given by $2r_{1}$ in the Global Plucker Formula $r_{k-1}-2r_{k}+r_{k+1}=2g-2-\beta_{k}$, where $g$ is the genus of the curve, $\beta_{k}$ is the ramification index of the $k$-th associated curve and $r_{k}$ is the number of osculating $k$-planes to the curve meeting a generic $(n-k-1)$-plane in $\mathbb{P}^{n}$(see \cite{genero}, Chapter 2). Since $r_{-1}=0$, $r_{0}$ is the degree of the curve and $\beta_{0} = 0$ (regularity condition), we have $2r_{1} = 4g-4+4d_{1}d_{2}$. As the genus of a curve of degree $d_{1}d_{2}$ is $g=\frac{d_{1}d_{2}(d_{1}+d_{2}-4)}{2}+1$, we have $2r_{1} = 2d_{1}d_{2}(d_{1}+d_{2}-2)$ isotropic points in $\mathbb{P}^{3}$. Assuming that the curve does not have isotropic points at infinity, all those points are in $\mathbb{C}^{3}$. 
\end{rem}

We now compute the global number of points where the osculating plane of an algebraic space curve is isotropic. 

\begin{prop} \label{oscul}
Let $f, g \in \mathbb{C}[x,y,z]$ be polynomials of degree $d_{1}$ and $d_{2}$, respectively, and consider the algebraic curve $\Phi= (f,g)$ defined by the intersection of the surfaces $f=0$ and $g=0$. {Assume that the projective closure of $\Phi$ is nonsigular at every point at infinity and intersects the plane at infinity  transversely. If the limiting osculating plane at every point at infinity is non-isotropic and the number $\operatorname{Osc}_{\Phi}$ of points of $\Phi$ with isotropic osculating planes is finite,  then $\operatorname{Osc}_{\Phi}= 6d_{1}d_{2}(d_{1}+d_{2}-3)$.}
\end{prop}

\begin{proof}
Let $\gamma: D \rightarrow \mathbb{C}^{3}, \gamma(t) = (x(t),y(t),z(t))$ be the parametrisation of the algebraic curve $\Phi$, where $D$ is an open subset of $\mathbb{C}$. 






Since $\gamma' = \nabla f \times \nabla g$, we have $\gamma''(t) = \frac{d}{dt}\nabla f(\gamma(t)) \times \nabla g + \nabla f \times \frac{d}{dt}\nabla g(\gamma(t))$, it follows from Lagrange's Formula that 
 $\gamma'(t) \times \gamma''(t)$ is parallel to
    \begin{align}\label{prod}
    q :=-[\mbox{det }(\nabla f, \nabla g, H(f) \cdot \gamma' )]\nabla g + [\mbox{ det }(\nabla f, \nabla g, H(g) \cdot \gamma')]\nabla f,
\end{align}
where $H(f)$ denotes the Hessian matrix of a polynomial $f$.

Hence, $Q:= \langle q, q \rangle$ is a polynomial of degree $ 6d_{1}+6d_{2}-14$.

Consider the homogenizations of $f$ and $g$ given by
\begin{align*}
    F &= f_{d_{1}} + wf_{d_{1-1}} + \dots + w^{d_{1}}f_{0}\\
    G &= g_{d_{2}} + wg_{d_{2-1}} + \dots + w^{d_{2}}g_{0},
\end{align*}
and still denote by $q$ the homogenization of the polynomial  in \eqref{prod}.
It follows by Bézout's theorem that the system
\begin{equation}\label{system}
\begin{cases}
    F&=0\\
    G&=0\\
   Q&=0
\end{cases}
\end{equation}
has $d_{1}d_{2}(6d_{1}+6d_{2}-14)$ solutions in $\mathbb{CP}^{3}$, corresponding to $d_{1}d_{2}(6d_{1}+6d_{2}-14)$ points in $\mathbb{CP}^{3}$ at which the osculating plane is isotropic. 

Let $(\alpha_{1}: \alpha_{2}: \alpha_{3}:0) \in \mathbb{CP}^{3}$ be a solution of the above system and assume $\alpha_{1} \neq 0$. By Euler's identity, we have
 \begin{align*}
 ({f_{d_{1}}}_{z}{g_{d_{2}}}_{x} - {f_{d_{1}}}_{x}{ g_{d_{2}}}_{z})  &=\frac{\alpha_{2}}{\alpha_{1}}({ f_{d_{1}}}_{ y}{g_{d_{2}}}_{z} - {f_{d_{1}}}_{z}{g_{d_{2}}}_{y});\\
({f_{d_{1}}}_{x}{g_{d_{2}}}_{y} - {f_{d_{1}}}_{y}{g_{d_{2}}}_{x}) &= \frac{\alpha_{3}}{\alpha_{1}}({f_{d_{1}}}_{y}{g_{d_{2}}}_{z} - {f_{d_{1}}}_{z}{g_{d_{2}}}_{y}).
 \end{align*}

At $\alpha$, we have
\begin{align*}
   Q= ({f_{d_{1}}}_{y}{g_{d_{2}}}_{z} - {f_{d_{1}}}_{z}{g_{d_{2}}}_{y})^2 ||(M_{1}+\frac{\alpha_{2}}{\alpha_{1}}M_{2} + \frac{\alpha_{3}}{\alpha_{1}}M_{3}) \nabla f_{d_{1}} - (L_{1}+\frac{\alpha_{2}}{\alpha_{1}}L_{2} + \frac{\alpha_{3}}{\alpha_{1}}L_{3})\nabla g_{d_{2}}||^{2},
\end{align*}
where $||v||^{2}= \langle v,v\rangle$ for $v \in \mathbb{C}^{3}$ and $L_{i}$ and $M_{i}$ are the i-th components of $H(f) \cdot \gamma'$ and $H(g) \cdot \gamma'$ at $\alpha$, respectively, that is,
\begin{align*}
    L_{1} &= ({f_{d_{1}}}_{xx}x'+{f_{d_{1}}}_{xy}y'+{f_{d_{1}}}_{xz}z')(\alpha);\\
    L_{2} &= ({f_{d_{1}}}_{xy}x'+{f_{d_{1}}}_{yy}y'+{f_{d_{1}}}_{yz}z')(\alpha);\\
     L_{3} &= ({f_{d_{1}}}_{zx}x'+{f_{d_{1}}}_{zy}y'+{f_{d_{1}}}_{zz}z')(\alpha)
\end{align*}
and 
\begin{align*}
    M_{1} &= {g_{d_{2}}}_{xx}x'+{g_{d_{2}}}_{xy}y'+{g_{d_{2}}}_{xz}z'(\alpha);\\
    M_{2} &= {g_{d_{2}}}_{xy}x'+{g_{d_{2}}}_{yy}y'+{g_{d_{2}}}_{yz}z'(\alpha);\\
     M_{3} &= {g_{d_{2}}}_{zx}x'+{g_{d_{2}}}_{zy}y'+{g_{d_{2}}}_{zz}z'(\alpha).
\end{align*}
It follows by Euler's identity that 
\begin{align*}
    {f_{d_{1}}}_{xx}(\alpha)+\frac{\alpha_{2}}{\alpha_{1}}{f_{d_{1}}}_{xy}(\alpha)+\frac{\alpha_{3}}{\alpha_{1}}{f_{d_{1}}}_{zx}(\alpha) &= \frac{(d_{1}-1)}{\alpha_{1}}{f_{d_{1}}}_{x}(\alpha);\\
    {f_{d_{1}}}_{xy}(\alpha)+\frac{\alpha_{2}}{\alpha_{1}}{f_{d_{1}}}_{yy}(\alpha)+\frac{\alpha_{3}}{\alpha_{1}}{f_{d_{1}}}_{zy}(\alpha) &= \frac{(d_{1}-1)}{\alpha_{1}}{f_{d_{1}}}_{y}(\alpha);\\
    {f_{d_{1}}}_{zx}(\alpha)+\frac{\alpha_{2}}{\alpha_{1}}{f_{d_{1}}}_{zy}(\alpha)+\frac{\alpha_{3}}{\alpha_{1}}{f_{d_{1}}}_{zz}(\alpha) &= \frac{(d_{1}-1)}{\alpha_{1}}{f_{d_{1}}}_{z}(\alpha).
\end{align*}
Substituting the three equations above in the expressions of $L_{1}, L_{2}$ and $L_{3}$ we obtain $L_{1}+\frac{\alpha_{2}}{\alpha_{1}}L_{2}+\frac{\alpha_{3}}{\alpha_{1}}L_{3} = 0$. Analogously, $M_{1}+\frac{\alpha_{2}}{\alpha_{1}}M_{2}+\frac{\alpha_{3}}{\alpha_{1}}M_{3} = 0$. Therefore, $Q(\alpha) = 0$. Consequently, the osculating plane is isotropic at every point at infinity. Since the curve intersects $w=0$ at $d_{1}d_{2}$ points, there are $d_{1}d_{2}$ points at infinity where the osculating plane of the curve is isotropic.  We compute the intersection multiplicity of the projective closure of the curve $\Phi$ and $Q=0$ at each of these points. 

Since the projective closure of $\Phi$ is a regular curve, we parametrise it as a graph of a function $\Gamma$ with $\Gamma(0) = \alpha$ and with a simple root at $w$. In the affine coordinate $x=1$, the curve $\Gamma$ represents a curve with a simple pole at infinity. Therefore, it can be parametrised by 
\begin{align*}
\gamma(s) = \frac{\alpha}{s} + a_{0}+a_{1}s+a_{2}s^{2}+ O(3).
\end{align*}

It follows by the homogenization of the polynomial $Q$ that $Q(x:y:z:w) = w^{\deg Q}Q(\frac{x}{w}, \frac{y}{w}, \frac{z}{w})$, where $\deg Q$ is the degree of the polynomial $Q$. The intersection multiplicity of the projective closure of $\Phi$ and $Q=0$ at the point $\alpha$ is then given by 
\begin{align}
    \deg Q - \mbox{order of the pole of } (Q \circ \gamma)(s) \mbox{ at }s=0.
\end{align}

We now look for the order of the pole of $Q \circ \gamma$ at the origin.
 

Since $\Phi(\gamma(s(t)))=0$, we have
\begin{align*}
    d\Phi(\gamma(s(t))) \cdot \gamma_{s}\cdot s' = 0.
\end{align*}
Therefore, the vector $\gamma_{s} \cdot s'$ is parallel to $(\nabla f \times \nabla g)(\gamma(s))$. 
Moreover, differentiating $(\nabla f \times \nabla g)(\gamma(s))$ with respect to $t$ yields
\begin{align} 
\mbox{H}({f}) \cdot \gamma' \cdot s' \times \nabla g+ \nabla f \times \mbox{H}({g})\cdot \gamma' \cdot s' = \gamma_{ss} \cdot s'^{2} + s'' \cdot \gamma_{s}.
\end{align}
Therefore,
\begin{align} \label{produtovet}
 (\gamma'\times \gamma'')(t)= {s'}^{3}(\gamma_{s}\times \gamma_{ss})(s) 
\end{align}
and 
\begin{align} \label{Q}
\langle \gamma'(t) \times \gamma''(t), \gamma'(t) \times \gamma''(t) \rangle = s'^{6} \langle \gamma_{s}(s) \times \gamma_{ss}(s), \gamma_{s}(s) \times \gamma_{ss}(s) \rangle.
\end{align}

Since $\gamma_{s} = \frac{-\alpha}{s^{2}}+ a_{1}+ 2a_{2}s+ O(s^{2})$ and $\gamma_{ss} = \frac{2 \alpha}{s^{3}} + 2a_{2} + O(s)$, we have
\begin{align*}
 \gamma_{s}\times \gamma_{ss} &= \left(\frac{-\alpha}{s^2}+ a_{1}+2a_{2}s+  O(s^{2})\right) \times \left(\frac{2\alpha}{s^{3}}+2a_{2}+O(s)\right).
 \end{align*}
 Hence, the limiting osculating plane at $\alpha$, namely, the plane through  $\alpha$  orthogonal to the limit of the direction of $\gamma_{s} \times \gamma_{ss}$ as $s \to 0$ is given by $a_{1} \times \alpha$ and 
 \begin{align}\label{result}
 \langle \gamma_{s} \times \gamma_{ss}, \gamma_{s} \times \gamma_{ss} \rangle = \frac{4}{s^{6}} \langle a_{1} \times \alpha, a_{1} \times \alpha \rangle + O(s^{-5}).
\end{align}

By \eqref{result}, $\langle \gamma_{s}\times \gamma_{ss}, \gamma_{s}\times \gamma_{ss} \rangle $ has a pole of order six, provided $\langle a_{1} \times \alpha, a_{1} \times \alpha \rangle \neq 0$, that is, the limiting osculating plane of the curve at a point at infinity is non-isotropic. 

We now compute the order of the pole of $s'$. Since the degree of $\nabla F \times \nabla G$ is $d_{1}+d_{2}-2$ and $\nabla F \times \nabla G = w^{d_{1}+d_{2}-2} (\nabla f \times \nabla g)$, we have that $\nabla f \times \nabla g$ has a pole of order $d_{1}+d_{2}-2$. Moreover, from
\begin{align*}
    (\nabla f \times \nabla g) (\gamma(s)) &= s' \gamma'(s)\\
    &= s'(\frac{-\alpha}{s^{2}}+a_{1}+2a_{2}s+ O(s^{2})),
\end{align*}
we obtain that order of the pole of $s'$ is $d_{1}+d_{2}-4$. By hypothesis, the limiting osculating plane is not isotropic at $\alpha$. So, using \eqref{result} we obtain that the order of the pole of $Q \circ \gamma$ at the origin is  $6(d_{1}+d_{2}-4) + 6 = (6d_{1}+6d_{2}-18)$, provided the limiting osculating plane is not isotropic at infinity. The intersection multiplicity of the projective closure of $\Phi$ with $Q=0$ is then given by $(6d_{1}+6d_{2}-14) - (6d_{1}+6d_{2}-18) = 4$.

Therefore, counted with multiplicity, we have $4d_{1}d_{2}$ points
at infinity whose osculating plane is isotropic. So, the curve $\Phi$ has $d_{1}d_{2}(6d_{1}+6d_{2}-14)-4d_{1}d_{2} = 6d_{1}d_{2}(d_{1}+d_{2}-3)$ points on which the osculating plane is isotropic.
\end{proof}

\begin{rem}
The number of points where the osculating plane is isotropic can also be calculated using the Plucker Formula. We first observe that the isotropic osculating plane is the osculating plane of the curve which is tangent to the conic $\{x^2+y^2+z^2=0; w=0\}$. Therefore, the number of points in $\mathbb{P}^{3}$ where the osculating plane is isotropic corresponds to the number of points where the osculating plane is tangent to that conic. By the general Plucker formula, we have $2r_{2} = 6d_{1}d_{2}(d_{1}+d_{2}-3)$ points of the curve in $\mathbb{P}^{3}$ where the osculating plane is isotropic. This coincides with the number of points on the affine curve $\Phi$ at which the osculating plane is isotropic, provided that no such points occur at infinity.
\end{rem}

The next proposition shows that the global number of zero-torsion points of an algebraic space curve depends on both the osculating plane and the singularities of the curve at infinity.

\begin{prop}
Let $f, g \in \mathbb{C}[x,y,z]$ be polynomials with degree $d_{1}$ and $d_{2}$, respectively, and let $\Phi= (f,g)$ be the algebraic curve given by the intersection of the surfaces $f=0$ and $g=0$, with a  finite number of zero-torsion points. Assume that the projective closure of $\Phi$ intersects the plane at infinity $w=0$ transversely and it is a nonsigular curve at infinity. Then, $\Phi$ has $\operatorname{T}_{\Phi}=2d_{1}d_{2}(3d_{1} + 3d_{2} - 10)$ zero-torsion points in $\mathbb{C}^{3}$, provided the projective closure of $\Phi$ does not have a zero-torsion point at infinity. 
\end{prop}
\begin{proof}
Let $\gamma: D \rightarrow \mathbb{C}^3$, $\gamma(t) = (x(t), y(t), z(t))$ be the parametrisation of the curve $\Phi=(f,g)$, where $D$ is an open subset of $\mathbb{C}$. 

    Differentiating $\gamma'(t) = \nabla f(\gamma(t)) \times \nabla g(\gamma(t))$ with respect to $t$, we obtain
    \begin{align*}
        \gamma''(t)&= \frac{d}{dt}(\nabla f(\gamma(t))) \times \nabla g(\gamma(t))+\nabla f(\gamma(t)) \times \frac{d}{dt}(\nabla g(\gamma(t)))\\
        \gamma'''(t)&= \frac{d^2}{dt^{2}}(\nabla f(\gamma(t)) \times \frac{d}{dt}\nabla g +2\frac{d}{dt}\nabla f(\gamma(t)) \times \nabla g(\gamma(t))+\nabla f(\gamma(t)) \times \frac{d^2}{dt^2} \nabla g(\gamma(t)).
    \end{align*}
    
Since $\gamma'(t) \times \gamma''(t) = - \beta_{1}\nabla g + \beta_{2}\nabla f$, where $\beta_{1}:= det[\nabla f, \nabla g, H(f) \cdot  \gamma']$ and $\beta_{2}:= det[\nabla f, \nabla g,  H(g) \cdot \gamma']$, we have that $\langle \gamma'(t) \times \gamma''(t), \gamma'''(t)\rangle$ is given by
\begin{align}\label{deter}
h:= \langle \beta_{2}\nabla f, \frac{d^2}{dt^{2}}\nabla f \times \nabla g \rangle 
+2\langle -\beta_{1}\nabla g + \beta_{2} \nabla f, H(f) \cdot \gamma' \times H(g) \cdot \gamma' \rangle \notag \\
+\langle -\beta_{1} \nabla g, \nabla f \times \frac{d^2}{dt^2} \nabla g \rangle. 
\end{align}
The zero-torsion points of the curve are then the solutions of the system 
\begin{equation}\label{torsion}
\begin{cases}
    f&=0\\
    g&=0\\
h&=0.
\end{cases}
\end{equation}
The homogenisation of the polynomials $f$ and $g$ yields the homogenisation $H$ of the polynomial $h$. We then apply Bézout's theorem in order to determine the number of intersection points of these curves in $\mathbb{PC}^{3}$ and consider the solutions lying at infinity. 

Let $F(x:y:z:w)=0$ and $G(x:y:z:w)=0$ be the homogenizations of $f=0$ and $g=0$ given by
\begin{align*}
    F(x:y:z:w)&=f_{d_{1}}(x,y,z)+wf_{d_{1}-1}(x,y,z)+ \cdots+ w^{d_{1}}f_{0}\\
    G(x:y:z:w)&= g_{d_{2}}(x,y,z)+wg_{d_{2}-1}(x,y,z)+\cdots + w^{d_{2}}g_{0}.
\end{align*}
By Bézout's Theorem, there are $d_{1}d_{2}(6d_{1}+6d_{2}-15)$ zero-torsion points in $\mathbb{CP}^{3}$.

Let $\alpha:=(\alpha_{1}:\alpha_{2}:\alpha_{3}:0) \in \mathbb{PC}^{3}$ be a solution of the homogenised system \eqref{torsion}. 


Since ${f_{d_{1}}}_{x}(\alpha) = -\frac{\alpha_{2}}{\alpha_{1}} {f_{d_{1}}}_{y}(\alpha) - \frac{\alpha_{3}}{\alpha_{1}}{f_{d_{1}}}_{z}(\alpha)$ and ${g_{d_{2}}}_{x}(\alpha) = - \frac{\alpha_{2}}{\alpha_{1}}(\alpha){g_{d_{2}}}_{y} - \frac{\alpha_{3}}{\alpha_{1}}{g_{d_{2}}}_{z}(\alpha)$, we have at $\alpha$ 

\begin{align}\label{bet1}
    \beta_{2}\langle \nabla f_{d_{1}}, \frac{d^{2}}{dt^{2}} \nabla {f}_{{d}_{1}} \times \nabla g_{d_{2}} \rangle &= (A_{1}+\frac
    {\alpha_{2}}{\alpha_{1}}A_{2}+\frac{\alpha_{3}}{\alpha_{1}}A_{3})(g_{d_{2_y}}f_{d_{1_z}}-f_{d_{1_y}}g_{d_{2_z}})
    \beta_{2}
    \end{align}
    \begin{align}\label{bet2}
 \langle \nabla {f_{d_1}} \times \frac{d^{2}}{dt^{2}} \nabla {g_{d_2}},  -\beta_{1} \nabla {g}_{{d}_{2}}\rangle &= (-B_{1}-\frac{\alpha_{2}}{\alpha_{1}}B_{2}-\frac{\alpha_{3}}{\alpha_{1}}B_{3})\beta_{1}(g_{d_{2_y}}f_{d_{1_z}}-g_{d_{2_z}}f_{d_{1_y}}),
    \end{align}
with $\frac{{d}^{2}}{dt^{2}}\nabla f_{{d}_{1}} (\alpha)= (A_{1}, A_{2}, A_{3})$ and $\frac{{d}^{2}}{dt^{2}}\nabla g_{{d}_{2}}(\alpha) = (B_{1}, B_{2}, B_{3})$. Differentiating $xF_{x}+yF_{y}+zF_{z}+wF_{w} = d_{1}F$ and evaluating at $\alpha$ we obtain 
\begin{align*}
f_{d_{1_{xx}}}+\frac{\alpha_{2}}{\alpha_{1}}f_{d_{1_{xy}}}+\frac{\alpha_{3}}{\alpha_{1}}f_{d_{1_{xz}}}=\frac{(d_{1}-1)}{\alpha_{1}}f_{d_{1_x}}\\
f_{d_{1_xy}}+\frac{\alpha_{2}}{\alpha_{1}}f_{d_{1_{yy}}}+\frac{\alpha_{3}}{\alpha_{1}}f_{d_{1_{zy}}}= \frac{(d_{1}-1)}{\alpha_{1}}f_{d_{1_y}}\\
f_{d_{1_{xz}}}+\frac{\alpha_{2}}{\alpha_{1}}f_{d_{1_{yz}}}+\frac{\alpha_{3}}{\alpha_{1}}f_{d_{1_{zy}}}= \frac{(d_{1}-1)}{\alpha_{1}}f_{d_{1_z}}.
\end{align*}
Moreover,
\begin{align*}
    f_{d_{1_{xxx}}}+\frac{\alpha_{2}}{\alpha_{1}}f_{d_{1_{xxy}}}+\frac{\alpha_{3}}{\alpha_{1}}f_{d_{1_{xxz}}}=\frac{(d_{1}-2)}{\alpha_{1}}f_{d_{1_{xx}}}.\\
f_{d_{1_{xyx}}}+\frac{\alpha_{2}}{\alpha_{1}}f_{d_{1_{xyy}}}+\frac{\alpha_{3}}{\alpha_{1}}f_{d_{1_{xzy}}}= \frac{(d_{1}-2)}{\alpha_{1}}f_{d_{1_{xy}}}.\\
f_{d_{1_{xxz}}}+\frac{\alpha_{2}} {\alpha_{1}}f_{d_{1_{xyz}}}+\frac{\alpha_{3}}{\alpha_{1}}f_{d_{1_{xzz}}}= \frac{(d_{1}-2)}{\alpha_{1}}f_{d_{1_{xz}}}.
\end{align*}
Substituting those derivatives in $A_{1}, A_{2}$ and $A_{3}$ we obtain $A_{1}+\frac{\alpha_{2}}{\alpha_{1}}A_{2}+\frac{\alpha_{3}}{\alpha_{1}}A_{3}$ is given by
\begin{align} \label{A1}
\frac{1}{\alpha_{1}}[(d_{1}-1)\langle \nabla f_{d_{1}}, H(f_{d_{1}}) \cdot \gamma' \times \nabla g_{d_{2}}+ \nabla f_{d_{1}} \times H(g_{d_{2}}) \cdot \gamma' \rangle + (d_{1}-2)\langle \gamma', H(f_{d_{1}}) \cdot \gamma' \rangle ].
    \end{align}
Analogously, $B_{1}+\frac{\alpha_{2}}{\alpha_{1}}B_{2}+\frac{\alpha_{3}}{\alpha_{1}}B_{3}$ is given by
    \begin{align}\label{B1}
\frac{1}{\alpha_{1}}[(d_{2}-1)\langle \nabla g_{d_{2}}, H(f_{d_{1}}) \cdot \gamma' \times \nabla g_{d_{2}}+ \nabla f_{d_{1}} \times H(g_{d_{2}}) \cdot \gamma'\rangle + (d_{2}-2) \langle \gamma', H(g_{d_{2}}) \cdot \gamma'\rangle].
\end{align}
As shown in the proof of Proposition \ref{oscul}, $L_{1} + \frac{\alpha_{2}}{\alpha_{1}}L_{2}+\frac{\alpha_{3}}{\alpha_{1}}L_{3} = M_{1} + \frac{\alpha_{2}}{\alpha_{1}}M_{2}+\frac{\alpha_{3}}{\alpha_{1}}M_{3}=0$, where $L_{i}$ and $M_{i}$ are i-th components of $H(f) \cdot \gamma'$ and $H(g) \cdot \gamma'$, respectively. Therefore, by \eqref{bet1} and \eqref{bet2} we have at $\alpha$: 


\begin{align}
     -\langle \nabla f_{d_{1}}, H(f_{d_{1}}) \cdot \gamma' \times \nabla g_{d_{2}}+ \nabla f_{d_{1}} \times H(g_{d_{2}}) \cdot \gamma' \rangle &= 0. \label{beta1}\\
- \langle \nabla g_{d_{2}}, H(f_{d_{1}}) \cdot \gamma' \times \nabla g_{d_{2}}+ \nabla f_{d_{1}} \times H(g_{d_{2}}) \cdot \gamma' \rangle &=0. \label{beta2}
\end{align}
Since $f_{d_{1_x}}= -\frac{\alpha_{2}}{\alpha_{1}}f_{d_{1_y}}-\frac{\alpha_{3}}{\alpha_{1}}f_{d_{1_z}}$ and $g_{d_{2_x}}= -\frac{\alpha_{2}}{\alpha_{1}}g_{d_{2_y}}-\frac{\alpha_{3}}{\alpha_{1}}g_{d_{2_z}}$, we have
\begin{align*}
    \langle \gamma', H(f_{d_{1}})\cdot \gamma'\rangle&=0.\\
     \langle \gamma', H(g_{d_{2}})\cdot \gamma'\rangle&=0.
\end{align*}
Thus, by \eqref{bet1}, \eqref{bet2}, \eqref{A1} and \eqref{B1} we have
\begin{align}
\beta_{2} \langle \nabla f_{d_{1}}, \frac{d^{2}}{dt^{2}} \nabla f_{d_{1}} \nabla g_{d_{2}} \rangle &= 0. \label{pfb}\\
\langle \nabla f_{d_{1}} \times \frac{d^{2}}{dt^{2}} \nabla g_{d_{2}}, -\beta_{1}\nabla g_{d_{2}} \rangle &= 0. \label{sfb}
\end{align}

It follows from $\eqref{beta1}$, $\eqref{beta2},\eqref{pfb}$ and $\eqref{sfb}$
that $\alpha$ is a solution of the homogenization of $h=0$, with $h$ as in \eqref{deter}. Hence, every point of the curve at infinity is a zero-torsion point. Since the curve intersects the plane $w=0$ at $d_{1}d_{2}$ points we have $d_{1}d_{2}$ zero-torsion points in $\mathbb{CP}^{3}$. We now compute the intersection multiplicity of the projective closure of the curve $\Phi$ and the homogenization of $h=0$ at a point at infinity, in an analogous way to that of Proposition \ref{oscul}.

As the projective closure of $\Phi$ is a regular curve at $\alpha$, we parametrised it as a graph of a function $\psi$ with a simple root at the coordinate $w=0$ and with $\psi(0) = \alpha$. Denoting by $\gamma$ its representation in the affine coordinate $x=1$, we have
\begin{align}\label{gamma}
    \gamma(s) = \frac{\alpha}{s}+ b_{0}+b_{1}s+ b_{2}s^{2}+b_{3}s^{3}+O(s^{4}).
\end{align}

The homogenisation of the polynomial $h$ yields $H(x:y:z:w) = w^{\deg H}h(\frac{x}{w}, \frac{y}{w}, \frac{z}{w})$, where $\deg H$ is the degree of the polynomial $H$. The intersection multiplicity of the projective closure of $\Phi$ with $H=0$ at the point $\alpha$, is then given by
\begin{align*}
    \deg H - \mbox{ order of the pole of } H \circ \gamma \mbox{ at }s = 0.
\end{align*}

We now compute the order of the pole of $H \circ \gamma$ at $s=0$.

Differentiating
\begin{align*}
    \gamma''(s(t)) = H_{f} \cdot \gamma' \cdot s' \times \nabla g + \nabla f \times H_{g} \cdot \gamma' \cdot s' = \gamma_{s} \cdot s'' +  s'^{2} \cdot \gamma_{ss}.
\end{align*}
with respect to $t$, yields
\begin{align}\label{gamma3}
    \gamma^{(3)}(s(t)) = s^{(3)} \cdot \gamma_{s} + 3s' \cdot s''\cdot \gamma_{ss} + s'^{3} \cdot \gamma_{sss}.
\end{align}

Since $(\gamma'\times \gamma'')(t) = s'^{3}(\gamma_{s} \times \gamma_{ss})$, it follows by \eqref{gamma3} that
\begin{align*}
    \langle \gamma'(t) \times \gamma''(t), \gamma'''(t) \rangle &= \langle s'^{3}(\gamma_{s} \times \gamma_{ss}), \gamma_{sss} \rangle\\
    &= s'^{6}\langle \gamma_{s} \times \gamma_{ss}, \gamma_{sss} \rangle.
\end{align*}

From \eqref{gamma}, we have 
\begin{align*}
    \gamma_{s}&= - \frac{\alpha}{s^{2}}+b_{1}+O(s).\\
    \gamma_{ss}&= \frac{2\alpha}{s^{3}}+b_{2}+6b_{3}s+O(s^{2}).\\
    \gamma_{sss}&= \frac{-6 \alpha}{s^{4}}+6b_{3}+24b_{4}s+O(s^{2}).
\end{align*}
Thus, $\langle \gamma_{s} \times \gamma_{ss}, \gamma_{sss} \rangle$ is given by
\begin{align*}
 \frac{-6}{s^{4}} \det(\beta', \beta'', \alpha) + \frac{2}{s^{3}}\det(\beta', \alpha, \beta''')-\frac{1}{s^{2}}\det(\alpha, \beta'', \beta''')+\det(\beta', \beta'', \beta''').
\end{align*}
Therefore, $\langle \gamma_{s} \times \gamma_{ss}, \gamma_{sss} \rangle$ has a pole of order four, provided that $ \det(\beta', \beta'', \alpha)= 2\det(\alpha, b_{1},b_{2}) \neq 0$, that is, $\alpha$ is not a zero-torsion point of the curve. As this condition holds by hypothesis, it follows that the order of the pole of $H \circ \gamma=s'^{6}\langle \gamma_{s} \times \gamma_{ss}, \gamma_{sss} \rangle$ is $6(d_{1}+d_{2}-4)+4= 6d_{1}+6d_{2}-20$, since the order of the pole of $s'$  is $d_{1}+d_{2}-4$, as shown in the proof of Proposition \ref{oscul}. The intersection multiplicity of the projective closure of $\Phi$ with $H = 0$ is then given by $(6d_{1}+6d_{2}-15)-(6d_{1}+6d_{2}-20)= 5$.  This yields that the number of zero-torsion points is $d_{1}d_{2}(6d_{1}+6d_{2}-15) - 5d_{1}d_{2} = 2d_{1}d_{2}(3d_{1}+3d_{2}-10)$, provided the curve does not have a zero-torsion point at infinity. 
\end{proof}
\begin{rem}
Since the zero-torsion points are the points where the osculating plane has contact at least four with the curve, the number of zero-torsion points is given by $\beta_{2}$ in the general Plucker formula. It follows from the fact that the curve is regular and with a non-degenerate contact with its tangent line that $\beta_{0}=\beta_{1}=0$. Hence, by the Plucker formula, we have $\beta_{2} = 2d_{1}d_{2}(3d_{1}+3d_{2}-10)$, which coincides with the number of zero-torsion points of the affine curve $\Phi$, provided that no zero-torsion points occur at infinity.
\end{rem}

\begin{cor}
Let $f,g \in \mathbb{C}[x,y,z]$ be two polynomials with degree $d_{1}$ and $d_{2}$, respectively and let $\Phi = (f,g)$ be the algebraic curve in $\mathbb{C}^{3}$ given by the intersection of the surfaces $f=0$ and $g=0$. Assume that no point at infinity is an isotropic point, a zero-torsion point, or a point with an isotropic osculating plane. If the projective closure of $\Phi$ is nonsingular at infinity and $\operatorname{Iso}_{\Phi},\operatorname{Osc}_{\Phi}$ and $\operatorname{T}_{\Phi}$ are finite, then, 
\begin{align*}
    2 \operatorname{T}_{\Phi}+3\operatorname{ Iso}_{\Phi}-3 \operatorname{ Osc}_{\Phi}=2d_{1}d_{2}.
\end{align*}
\end{cor}

\subsection{Local geometric invariants for holomorphic curves}
We compute the number of isotropic points, zero-torsion points and points with isotropic osculating plane that appear under a deformation of a holomorphic space curve with a singular point at the origin.

\begin{defn}
For a germ of a nonzero holomorphic space curve with a singular point at the origin and parametrised by $\gamma:(\mathbb{C},0) \rightarrow (\mathbb{C}^{3},0)$, $\gamma(t) = (a_{k}t^{k} + \mbox{h.o.t.}, a_{m}t^{m} + \mbox{h.o.t.}, a_{n}t^{n} + \mbox{h.o.t.})$ with $k<m,n$ and $a_{k} \neq 0$, the limiting tangent direction of $\gamma$ at the origin is $(a_{k},0,0).$ 
\end{defn}

\begin{defn}
Let $\gamma:(\mathbb{C},0)\to(\mathbb{C}^{3},0)$ be a germ of a nonzero holomorphic space curve with a singular point at the origin. The limiting tangent vector of $\gamma$ at the origin, denoted by $\gamma'(0)$, is the limit of the  direction of $\gamma'(t)$ as $t \to 0$. The line through $\gamma(0)$ parallel to $\gamma'(0)$ is called the limiting tangent line of $\gamma$ at the origin. 
Moreover, we define the vector $\gamma'(0) \times \gamma''(0)$ as the limit direction of ${\gamma}'\times {\gamma}''$
at the origin. The limiting osculating plane of $\gamma$ at the origin is the plane through $\gamma(0)$ orthogonal to $\gamma'(0)\times {\gamma}''(0)$ and the limiting rectifying plane is the plane through $\gamma(0)$ and generated by $\gamma'(0)$ and $\gamma'(0) \times \gamma''(0)$.
\end{defn}

We first characterize the parametrisation of the curve in terms of its order of contact with its limiting rectifying plane and its limiting osculating plane at the origin, as well as its multiplicity.
\begin{theo}\label{carac}
Let  $\gamma: (\mathbb{C},0) \rightarrow (\mathbb{C}^{3},0)$ be a germ of a nonzero holomorphic space curve parametrised by $\gamma(t) = (t^{k}+at^{l_{{0}}} + \mbox{h.o.t.}, bt^{m_{0}} + \mbox{h.o.t.}, ct^{n_{0}} + \mbox{h.o.t.})$ with $l_{0}, m_{0}, n_{0} > k$ and $m_{0} \neq n_{0}$. The order  of contact of $\gamma$ with its limiting osculating plane at the origin is $\max\{n_{0}, m_{0}\}$ whereas its order of contact with the limiting rectifying plane at the origin is $\min\{n_{0}, m_{0}\}$.

\end{theo}

\begin{proof}
We suppose without loss of generality that $n_{0}> m_{0}$. The remaining case follows by an analogous argument.


Note that $  \gamma'(t) \times \gamma''(t)$ is given by
\begin{align*}
(At^{m_{0}+n_{0}-3} + \mbox{h.o.t.},
    Bt^{k+n_{0}-3}-Ct^{l_{0}+n_{0}-3} + \mbox{h.o.t.},Dt^{k+m_{0}-3}+Et^{l_{0}+m_{0}-3} + \mbox{h.o.t.}),
\end{align*}
with
\begin{align*}
    A&=bcm_{0}n_{0}(n_{0}-m_{0});
    &B&=-cn_{0}(n_{0}-k)k;\\
    C&=acn_{0}(n_{0}-l_{0})l_{0};
  &D&=b{m_{0}}(m_{0}-k)k;
 &E&=abm_{0}(m_{0}-l_{0})l_{0}.
\end{align*}
Since $n_{0}>m_{0}$ and $l_{0}, m_{0}, n_{0}>k$, the orthogonal vector to the osculating plane of the curve at the origin is then parallel to $(0,0,1)$.



The result follows by observing that the contact of the curve with this plane at the origin is measured by the singularities of the height function
\begin{align*}
 H: (\mathbb{C},0) &\rightarrow (\mathbb{C},0)\\
    t &\mapsto H(t) = \langle \gamma(t), \gamma'(0) \times \gamma''(0) \rangle=\langle \gamma(t) , (0,0,1) \rangle=t^{n_{0}} + \mbox{h.o.t.}
\end{align*}
at the origin. Hence, the order of contact of the curve and its osculating plane at the origin is $n_{0}= \max\{m_{0},n_{0}\}$.

We now compute the order of contact between the curve and its limiting rectifying plane at the origin. Since
\begin{align*}
    \gamma'(t)= (kt^{k-1}+al_{0}t^{l_{0}-1} + \mbox{h.o.t.}, bm_{0}t^{m_{0}-1} + \mbox{h.o.t.}, cn_{0}t^{n_{0}-1} + \mbox{h.o.t.}),
\end{align*}
the vector $\gamma'(0)$ is parallel to $(1,0,0)$. Since the limiting rectifying plane at the origin is generated by the vectors $(1,0,0)$ and $(0,0,1)$, it follows that its orthogonal vector at that point is parallel to $(0,1,0)$. The contact of the curve with its limiting rectifying plane at the origin is then given by the singularities of the height function
\begin{align*}
H: (\mathbb{C},0) &\rightarrow (\mathbb{C},0)\\
t &\mapsto H(t)= \langle \gamma(t), (0,1,0)\rangle = t^{m_{0}} + \mbox{h.o.t.}.
\end{align*}
at the origin. Thus, the order of contact between the curve and its limiting rectifying plane at the origin is  $m_{0} = \min\{n_{0}, m_{0}\}$.

\end{proof}

Using the same parametrisation as above, the number of isotropic points, zero-torsion points and points with isotropic osculating plane that appear after a deformation of a germ of a curve $\gamma: (\mathbb{C},0) \rightarrow (\mathbb{C}^{3},0)$ are characterized, respectively, by the orders of contact at the origin between the curve and its rectifying plane and between the curve and its osculating plane, as follows.

\begin{theo}
Let $\gamma: (\mathbb{C},0) \rightarrow (\mathbb{C}^{3},0)$ be a germ of a nonzero holomorphic curve with multiplicity $k$ and a singularity at the origin. Suppose that the limiting tangent line {and the limiting osculating plane} of the curve at the origin are non-isotropic and the orders of contact of $\gamma$ with its limiting rectifying plane and with its limiting osculating plane at the origin are $r$ and $s$, respectively. Then,
\begin{align*}
\operatorname{Iso}_{\gamma} &= 2k-2\\
\operatorname{Osc}_{\gamma} &= 2k-6+2r\\
    \operatorname{T}_{\gamma}&=-6+k+r+s.
\end{align*}
\end{theo}
\begin{proof}
We parametrise the germ of the curve by $\gamma: (\mathbb{C},0) \rightarrow (\mathbb{C}^{3},0), \gamma(t) = (t^{k}+at^{l_{0}} + \mbox{h.o.t.}, bt^{m_{0}} + \mbox{h.o.t.}, ct^{n_{0}} + \mbox{h.o.t.})$, with $l_{0}, m_{0}, n_{0}>k$ {and $m_{0} \neq n_{0}$}. Note that
\begin{align*}
    \operatorname{Iso}_{\gamma}&=\operatorname{ord}(\langle \gamma'(t), \gamma'(t) \rangle) \\
    &= \operatorname{ord}(k^{2} t^{-2 + 2 k}+ 2 akl_{0} t^{-2 + k + l0} + a^{2} l_{0}^{2} t^{-2 + 2l_{0}} +
  b^{2}m_{0}^2t^{-2 + 2 m_{0}} + c^{2}n_{0}^{2}t^{-2+2n_{0}})\\
  &= 2k-2.
\end{align*}

Moreover, since
\begin{align*}
    \langle \gamma'(t) \times \gamma''(t), \gamma'(t) \times \gamma''(t) \rangle &= At^{-6+2k+2 m_{0}}+Bt^{-6+2 k+2 n_{0}} + \mbox{h.o.t.},
    \end{align*}
with 
  \begin{align*}
  A&=b^{2}k^{2}m_{0}^{2}(k-m_{0})^{2}\\
  B&=c^{2}k^{2}n_{0}^{2}(k-n_{0})^{2}
  \end{align*}  
we obtain 
\begin{align*}
\operatorname{Osc}_{\gamma}&=\operatorname{ord}(\langle \gamma'(t) \times \gamma''(t), \gamma'(t)\times \gamma''(t)\rangle)
= 2k-6+2\min\{m_{0}, n_{0}\} = 2k-6+2r,
\end{align*}
where the last equality follows from Theorem \ref{carac}.

Finally, we compute the number of zero-torsion points concentrated at the origin. We observe that the term of the lowest order of $\det(\gamma'(t), \gamma''(t), \gamma'''(t))$ is 
\begin{align*}
=bckm_{0}n_{0}(-k^{2}m_{0}+km_{0}^2+k^{2}n_{0}-m_{0}^{2}n_{0}-kn_{0}^{3}+m_{0}n_{0}^{2})t^{-6 + k + m_{0} + n_{0}}.
\end{align*}
Therefore, 
\begin{align*}
    \operatorname{T}_{\gamma}= \operatorname{ord}(\det(\gamma'(t), \gamma''(t), \gamma'''(t)) = -6+k+m_{0}+n_{0}.
\end{align*}

Since $\min\{m_{0},n_{0}\}$ and $\max\{m_{0},n_{0}\}$ are the orders of contact of the curve with its rectifying plane and the osculating plane, respectively, it follows that 
\begin{align*}
    \operatorname{T}_{\gamma} = -6+k+r+s.
\end{align*}
\end{proof}

We now consider germs of singular curves whose limiting tangent direction at the origin is isotropic. In this case, the curve is tangent to the isotropic plane orthogonal to its limiting tangent vector at the origin. This plane is called the isotropic plane of the curve at the origin.

\begin{theo}
    Consider a germ of a singular curve with multiplicity $k$ and suppose that the limiting tangent direction of the curve at the origin is an isotropic direction. If $\gamma: (\mathbb{C},0) \rightarrow (\mathbb{C}^{3},0), \gamma(t) =  (t^{k}+at^{l_{0}} + \mbox{h.o.t.}, it^{k}+bt^{m_{0}} + \mbox{h.o.t.}, ct^{n_{0}} + \mbox{h.o.t.})$ is the parametrisation of the curve,  with $l_{0},m_{0}, n_{0}>k$ and $m_{0} \neq n_{0}$, the order of contact of the curve with its isotropic plane at the origin is $\min\{l_{0}, m_{0}\}$. Moreover, the curve has $\max\{n_{0},l_{0}, m_{0}\}$-order of contact with its osculating plane at the origin, whereas the order of contact between the curve and the isotropic cone at the origin is $\min\{l_{0}+k, m_{0}+k, 2n_{0}\}$.
\end{theo}
\begin{proof}
We suppose without loss of generality that $n_{0}>l_{0}, m_{0}$. The other case follows analogously.

We first compute the contact of the curve with its isotropic plane at the origin. 

Since
\begin{align*}
    \gamma'(t)&= (kt^{k-1}+al_{0}t^{l_{0}-1} + \mbox{h.o.t.},ikt^{k-1}+bm_{0}t^{m_{0}-1} + \mbox{h.o.t.},cn_{0}t^{n_{0}-1} + \mbox{h.o.t.}),
\end{align*}
the limiting tangent vector at the origin is parallel to $(1,i,0)$.



Note that the height function on $\gamma$ along the vector $(1,i,0)$ is described by
\begin{align*}
    H: (\mathbb{C},0) &\rightarrow (\mathbb{C},0)\\
    &t \mapsto H(t) = \langle\gamma(t), (1,i,0) \rangle=at^{l_{0}}+ibt^{m_{0}} + \mbox{h.o.t.}
\end{align*}
Hence, the order of contact of $\gamma$ with its isotropic plane is given by $\min\{l_{0},m_{0}\}$.

We now compute the order of contact of the curve and its osculating plane at the origin. Since
\begin{align*}
  \gamma'(t) \times \gamma''(t) =(ickn_0(n_0-k)t^{k+n_0-3}
+
bc\,m_0n_0(n_0-m_0)t^{m_0+n_0-3}
+ \mathrm{h.o.t.},
\\[2mm]
-ckn_{0}(n_0-k)t^{k+n_0-3}
-
ac\,l_0n_0(n_0-l_0)t^{l_0+n_0-3}
+ \mathrm{h.o.t.},
\\[2mm]
-iakl_{0}(l_0-k)t^{k+l_0-3}
+
bkm_{0}(m_0-k)t^{k+m_0-3}
+
ab\,l_0m_0(m_0-l_0)t^{l_0+m_0-3}
+\mathrm{h.o.t.})
\end{align*}
and $n_{0}>m_{0}$, the vector $\gamma'(0) \times \gamma''(0)$ is parallel to $(0,0,1)$. Hence, the order of contact of $\gamma$ with its osculating plane is $n_{0}$, since the height function along $\gamma'(0) \times \gamma''(0)$ has order $n_{0}$. 

The order of contact between the curve and the isotropic cone at the origin follows from observing the order of
\begin{align*}
D_{0}: (\mathbb{C},0) &\rightarrow (\mathbb{C},0)\\
t &\mapsto D_{0}(t)= \langle \gamma(t), \gamma(t) \rangle,
\end{align*}
where $D_{0}$ is the distance-squared function with respect to a circle centered at the origin. 
\end{proof}

We now investigate the relation between zero-torsion points, isotropic points and points where the limiting osculating plane is isotropic in the case of a germ of a singular curve with isotropic limiting tangent direction at the origin. 

\begin{theo}
    Consider a germ of a singular space curve $\gamma:(\mathbb{C},0) \rightarrow (\mathbb{C}^{3},0)$ with multiplicity $k$ and assume that the limiting tangent direction at the origin is parallel to an isotropic vector. Suppose that the orders of contact with the isotropic plane, the limiting osculating plane and the isotropic cone at the origin are $r$, $s$ and $l$, respectively. Then,
\begin{align*}
    \operatorname{Iso}_{\gamma}&=l-2\\
    \operatorname{Osc}_{\gamma}&=-6+k+l+r\\
    \operatorname{T}_{\gamma} &= -6+k+r+s.
\end{align*}
\end{theo}

\begin{proof}
We parametrise the germ of the curve as $\gamma: (\mathbb{C},0) \rightarrow (\mathbb{C}^{3},0), \gamma(t) = (t^{k}+at^{l_{0}} + \mbox{h.o.t.}, it^{k}+bt^{m_{0}} + \mbox{h.o.t.}, ct^{n_{0}} + \mbox{h.o.t.})$, with $l_{0},m_{0},n_{0}>k$ and suppose, without loss of generality, that $n_{0}>l_{0},m_{0}$. The other case is analogous. 

Note that
\begin{align*}
\langle \gamma'(t), \gamma'(t)\rangle &= 2akl_{0}t^{-2+k+l_{0}} + a^{2}l_{0}^{2} t^{2l_{0} - 2} +2ibkm_{0} t^{-2+k+m_{0}}+ b^{2} m_{0}^{2} t^{-2+2m_{0}}\\&+ c^{2} n_{0}^{2} t^{-2+2 n_{0}} + \mbox{h.o.t.}
\end{align*}
Hence,
\begin{align*}
\operatorname{Iso}_{\gamma} &= \operatorname{ord}(\langle \gamma'(t), \gamma'(t)\rangle)\\
&= \min \{l_{0}+k, m_{0}+k, 2n_{0}\}-2\\
&= l-2
\end{align*}

Now we compute the number of points concentrated at the origin where the osculating plane is isotropic. Since 
\begin{align*}
    \langle \gamma'(t) \times \gamma''(t), \gamma'(t) \times \gamma''(t) \rangle &= At^{-6+2k+2l_{0}}+Bt^{-6+2k+l_{0}+m_{0}}+Ct^{-6+2k+2m_{0}}\\
&+Dt^{-6+k+l_{0}+2n_{0}}+Et^{-6+k+m_{0}+2n_{0}} + \mbox{h.o.t.},
    \end{align*}
with
\begin{align*}
A&= -a^{2}k^{2}l_{0}^{2}(k^{2}-2kl_{0}+l_{0}^{2});\\
B&= -2iabk^{2}l_{0}m_{0}(k^{2}-kl_{0}-km_{0}+l_{0}m_{0});\\
C&= b^{2}k^{2}m_{0}^{2}(k^{2}-2km_{0}+m_{0}^{2});\\
D&=2ac^{2}kl_{0}n_{0}^{2}(kl_{0}-kn_{0}-l_{0}n_{0}+n_{0}^{2});\\
E&=2ibc^{2}km_{0}n_{0}^{2}(km_{0}-kn_{0}-m_{0}n_{0}+n_{0}^{2}).
\end{align*}

Hence, $\operatorname{Osc}_{\gamma} = -6+k+\min\{l_{0}, m_{0}\}$$+\min\{2n_{0}, l_{0}+k, m_{0}+k\} = -6+k+l+r$.
 
Now we compute the number of zero-torsion points concentrated at the origin. We observe that
\begin{align*}
    \det(\gamma'(t),\gamma''(t),\gamma'''(t))&= At^{-6+k+l_{0}+n_{0}} +Bt^{-6+k+m_{0}+n_{0}} + \mbox{h.o.t.},
\end{align*}
with
\begin{align*}
A &=iackl_{0}n_{0}(k^{2}l_{0}-kl_{0}^{2}-k^{2}n_{0}+l_{0}^{2}n_{0}-l_{0}n_{0}^{2}+kn_0^{2})\\
B &= -bckm_{0}n_{0}(k^{2}m_{0}-km_{0}^{2}-k^{2}n_{0}+m_{0}^{2}n_{0}+kn_{0}^{2}-m_{0}n_{0}^{2}).
\end{align*}
 Therefore, $\operatorname{T}_{\gamma} = -6+k+r+s$.
\end{proof}
\begin{rem}
The relation between the isotropic points, zero-torsion points and those points where the osculating plane is isotropic remains unchanged if we consider $k=1$, that is, for regular curves with an isotropic point at the origin. 
\end{rem}

\textbf{Acknowledgement:} The authors thank Farid Tari for valuable conversations and suggestions.

\textbf{Funding}: A.D.F. was supported by São Paulo Research Foundation (FAPESP), process
number 2025/12551-9. J.J.N.B. and R.O.S. were supported by grant PID2025-172404NB-I00 funded by MCIN/AEI and by “ERDF A way of making Europe”.






\end{document}